\documentclass[11pt]{article}
\usepackage{caption}
\usepackage{url}
\usepackage{verbatim}
\usepackage{enumitem}
\usepackage{amsmath}
\usepackage[usenames,dvipsnames]{color}
\usepackage{amsthm}
\usepackage{amsfonts}
\usepackage{graphicx}
\usepackage{nicefrac}
\usepackage{url}
\usepackage{etoolbox}
\usepackage{cmbright}
\usepackage{algorithm}
\usepackage{algpseudocode}
\usepackage{graphicx, amssymb, subcaption,graphics}
\newtheorem{theorem}{Theorem}[section]
\newtheorem{lemma}[theorem]{Lemma}

\theoremstyle{definition}

\theoremstyle{remark}
\newtheorem{remark}[theorem]{Remark}
\numberwithin{equation}{section}
\def\ti{\tilde}
\begin{document}

\title{An Efficient Finite Difference Scheme for the 2D Sine-Gordon Equation}

\author{
Xiaorong Kang\thanks{School of Science, Southwest University of Science and Technology, Mianyang, Sichuan 621010,  China (kangxiaorong@swust.edu.cn)}
\and
Wenqiang Feng\thanks{Department of Mathematics, The University of Tennessee, Knoxville, TN 37996 (wfeng1@utk.edu)}	
\and
Kelong Cheng\thanks{School of Science, Southwest University of Science and Technology, Mianyang, Sichuan 621010,  China (zhengkelong@swust.edu.cn)}
\and
Chunxiang Guo\thanks{School of Business, Sichuan University, Chengdu, Sichuan 610064, China (guocx70@163.com)}
}
\maketitle
\numberwithin{equation}{section}

\begin{abstract}
We present an efficient second-order  finite difference scheme for solving the 2D sine-Gordon equation,
which can inherit the discrete energy conservation for the undamped model theoretically. Due to the
semi-implicit treatment for the nonlinear term, it leads to a sequence of nonlinear coupled equations.
We use a linear iteration algorithm, which can solve them  efficiently, and the contraction mapping
property is also proven. Based on truncation errors of the numerical scheme, the convergence analysis
in the discrete $l^2$-norm  is investigated in detail. Moreover, we carry out various numerical
simulations, such as verifications of the second order accuracy, tests of energy conservation and circular
ring solitons, to demonstrate the efficiency and the robustness of the proposed scheme. 

\textbf{keyword} 2D sine-Gordon equation, conservative, difference scheme, linear iteration, convergence.

\end{abstract}

\maketitle




\section{Introduction}
\label{intro}

In this paper, we consider the following 2D sine-Gordon equation,
\begin{eqnarray}
u_{tt}+\beta u_t-\alpha \Delta u=-\phi(x,y)\sin u+F(x,y,t),\quad (x,y)\in \Omega, \quad t\ge 0,\label{1.1}
\end{eqnarray}
with initial conditions
\begin{eqnarray}
u(x,y,0)=\varphi_1(x,y),\quad
u_t(x,y,0)=\varphi_2(x,y),\quad (x,y)\in \Omega,\label{1.2}
\end{eqnarray}
and the boundary condition
\begin{eqnarray}
u|_{\partial \Omega}=G(t),\quad t\ge 0, \label{1.3}
\end{eqnarray}
where $\Omega=[0,L]^2$. This equation has attracted much attention  due to the presence
of soliton solutions and has  a great deal of  applications in the propagation of fluxons in
Josephson junctions between two superconductors \cite{Perring1962},  the motion of a rigid pendulum attached to a stretched wire \cite{Whitham1999},
dislocations in crystals  and the stability of fluid motions. Nowadays, it has become one of  paradigms of the nonlinear
dynamical system to describe many different physical phenomena \cite{Xin2000}. In (\ref{1.1}),  $\phi(x, y)$ is a
nonnegative function with finite bound $\phi_0$ and may be interpreted as the Josephson current density,
while  $\varphi_1(x, y)$ and $\varphi_2(x, y)$
represent wave modes or kinks and velocity, respectively. In particular, when $\beta=0$,  (\ref{1.1}) reduces to the undamped sine-Gordon equation,
\begin{eqnarray}
u_{tt}-\alpha \Delta u=-\phi(x,y)\sin u+F(x,y,t).  \label{1.4}
\end{eqnarray}
If $F=0$ and $G$ is periodic or homogeneous, one of the main properties of the undamped sine-Gordon
equation (\ref{1.4}) has the conservation for the energy defined as follows\cite{Bratsos2006, Bratsos2007},
\begin{eqnarray}
{ E}(t)&=&\frac{1}{2}\int_\Omega \big [|u_t|^2+|\nabla   u|^2+2\phi (1-\cos u)\big ]dxdy \nonumber\\
&=&\frac{1}{2}\big (\|u_t\|^2_{L^2(\Omega)}+\|\nabla u\|^2_{L^2(\Omega)}\big )+\int_\Omega \phi (1-\cos u)dxdy,\label{1.5}
\end{eqnarray}
which is not valid for the  damped system (\ref{1.1}).

Recently, various analytical and numerical methods have been proposed for the numerical solution
of partial differential equations, for example, the integral transform\cite{Liang2017, Yang20161, Yang20162, Yang2017}
and traveling-wave technologies\cite{Yang20172, Yang20163}. Analytical solutions to the unperturbed sine-Gordon equation
with zero damping  have been obtained by Lambs methods \cite{Zag1979} and B\"{a}cklund
transformations. Many efforts have been attempted to develop numerical methods,
such as the finite difference method \cite{Bratsos2006, Bratsos2007, Cui2010},
the time-splitting pseudospectral and spectral method \cite{Asg2013},
the finite element method \cite{Arg1991}, the  mesh-free reproducing kernel particle Ritz method \cite{Cheng2012},
the local weak meshless method \cite{Deh2010}, the boundary element method \cite{Deh2008-1}, the  Differential
quadrature method \cite{ Jiw2012} and the radial basis functions method \cite{Deh2008-2} for the 2D sine-Gordon equation.
However, there exist few available error estimate results in the above-mentioned works for the 2D case.
The main reason is that the techniques used for 1D case can not be extended trivially to high dimensions
because of the difficulty in obtaining the a priori uniform estimate of the numerical solution.

Since the undamped sine-Gordon equation is a conservative system, it should be pointed out that a conservative numerical
scheme performs better  than a nonconservative one. The key is that it can preserve some invariant
properties of the differential equation and capture physical procedures with more details \cite{Li1995}.
Moreover, there has been growing interest in conservative numerical methods for solving
partial differential equations. For example, Klein-Gordon equation \cite{Stra1978}, the high frequency wave
phenomena \cite{Cheng2015, Hu2014,  Wang2015, Zheng2013}, the
phase field crystal model \cite{Wang2011}, and so forth. As for the sine-Gordon equation, to our knowledge, only a few results
considered this vital property. Although many verification results of the discrete energy for various numerical
methods are reported\cite{Asg2013, Bratsos2006, Bratsos2007}, the analysis at a theoretical level was hardly shown.

The main purpose of this paper is to present a second-order semi-implicit finite difference
scheme for  numerical solutions of the 2D sine-Gordon equation (\ref{1.1})-(\ref{1.3}). There are three
main features to this work. The first is that the proposed scheme can admit the discrete energy
conversation for the undamped case (\ref{1.4}) at a theoretical level, which has not yet been
reported in the existing literatures. The second feature of this work is the linear iteration
algorithm introduced\cite{Wang2015} to solve efficiently the nonlinear system at each time step
due to the implicit treatment of the nonlinear term. Meanwhile, a careful analysis shows a
contraction mapping property of this iteration under the given constrain for the time step.
Finally, we provided a detailed convergence analysis for the second-order scheme in the $l^2$-norm.

The remainder of the paper is organized as follows. In Section \ref{S2}, the second-order finite
difference scheme is proposed and the energy conservation property for the undamped system is proven.
The  linear iteration algorithm and the corresponding  theoretical analysis of the contraction mapping are given in
Section \ref{S3}.  Truncation errors and the convergence analysis  are discussed
in Section \ref{S4}. Some numerical simulation results  are given to demonstrate the efficiency of
the linear iteration solver and the convergence of the scheme  in Section \ref{S5}. Finally, some
conclusions are made in Section \ref{S6}.

\section{Numerical scheme and energy conservation }
\label{S2}

\subsection{Second order finite difference scheme}

Let $v=u_t$.  (\ref{1.1}) can be rewritten as,
\begin{eqnarray}
&&v_{t}+\beta u_t -\alpha \Delta u=-\phi(x,y)\sin u+F(x,y,t),\label{2.1.1} \\
&&v=u_t.  \label{2.1.2}
\end{eqnarray}
Then, for the given 2D domain $\Omega$, define the uniform numerical grid $(x_i,y_j)$
with $\Delta x=\Delta y= h$ for simplicity of presentation. Let $M_x=M_y=M$ and $M\cdot h=L$ such
that $x_i=ih,i=0,1,\cdots,M, y_j=jh,j=0,1,\cdots,M$.  For a fixed time $T$, let $\Delta t$ be the
step size for temporal direction, $t^n=n\Delta t $, $n=0,1,2,\cdots,N$, $N=[\frac{T}{\Delta t}]$, $u_{i,j}^n\approx u(x_i,y_j, t^n)$.
Denote $\Delta_h=D_{xx}+D_{yy}$ as the standard second
order difference operator with
\begin{eqnarray*}
&&D_x u=\frac{u_{i+1,j}-u_{i,j}}{h},\quad D_y u=\frac{u_{i,j+1}-u_{i,j}}{h},\\
&&D_{xx} u=\frac{u_{i+1,j}-2u_{i,j}+u_{i-1,j}}{h^2},\quad  D_{yy} u=\frac{u_{i,j+1}-2u_{i,j}+u_{i,j-1}}{h^2}.
\end{eqnarray*}

The second order finite difference scheme is presented at a point-wise level as follows,
\begin{eqnarray}
&&\frac{v^{n+1}-v^n}{\Delta t}+\beta \frac{u^{n+1}-u^{n}}{\Delta t}-\frac{\alpha}{2}\Delta_h (u^{n+1}+u^n)=\phi \frac{\cos( u^{n+1})-\cos (u^n)}{u^{n+1}-u^n}+F^{n+\frac{1}{2}},\label{2.1.3}\\
&&\frac{u^{n+1}-u^{n}}{\Delta t}=\frac{v^{n+1}+v^n}{2},   \label{2.1.4}
\end{eqnarray}
with  discrete initial conditions
\begin{eqnarray}
u^0_{i,j}=\varphi_1(x_i,y_j),\qquad v^0_{i,j}=\varphi_2(x_i,y_j),  \label{2.1.5}
\end{eqnarray}
and the boundary condition
\begin{eqnarray}
u^n_{i,j}|_{\partial \Omega}=G(x_i,y_j,t^n), \qquad (x_i,y_j)\in \partial \Omega,  \label{2.1.6}
\end{eqnarray}
where $F^{n+\frac{1}{2}}=F(x_i,y_j,t^{n+\frac{1}{2}})$.

Obviously, (\ref{2.1.4}) can be reformulated as
\begin{eqnarray}
v^{n+1}=\frac{2(u^{n+1}-u^n)}{\Delta t}-v^n. \label{2.1.7}
\end{eqnarray}
Substituting (\ref{2.1.7}) into (\ref{2.1.3}) yields that
\begin{eqnarray}
\frac{2u^{n+1}}{\Delta t^2}+\frac{\beta}{\Delta t}u^{n+1}-\frac{\alpha}{2}\Delta_h u^{n+1}=\phi \frac{\cos( u^{n+1})
-\cos (u^n)}{u^{n+1}-u^n}+\frac{\beta}{\Delta t}u^{n}+\frac{\alpha}{2}\Delta_h u^{n}+\kappa (u^n,v^n)+F^{n+\frac{1}{2}},\label{2.1.8}
\end{eqnarray}
where $\kappa (u^n,v^n)=\frac{\frac{2u^n}{\Delta t}+2v^n}{\Delta t}$. (\ref{2.1.8}) is
nonlinear and can be solved implicitly by a linear iteration algorithm introduced
in the next section. Following $u^{n+1}$ is solved, $v^{n+1}$ can be computed explicitly by (\ref{2.1.7}).

\begin{remark}\label{rm0}
The main idea to deal with the sine nonlinearity  was first introduced by W. A. Strauss and
L. V\'{a}zquez\cite{Stra1978} to compute  numerical solutions of a nonlinear Klein-Gordon equation
in which a polynomial nonlinear term is involved. In fact, this subtle technique can achieve perfect
numerical solutions for the conservative model and has been extensively studied for some nonlinear problems,
such as the Cahn-Hilliard type equation\cite{Chen2016, Cheng2016}.
\end{remark}

\subsection{Discrete energy conservation for the undamped equation}
As above mentioned,  the undamped sine-Gordon equation with certain boundary conditions admits the property
of energy conservation. Here, as the special case of the numerical scheme (\ref{2.1.3}) and (\ref{2.1.4}),
the difference scheme for the undamped equation  is  conservative for the discrete energy.

Letting $F(x,y,t)=0$ and $G(t)=0$, from (\ref{1.4}), we have
\begin{eqnarray}
u_{tt}-\alpha \Delta u=-\phi(x,y)\sin u,  \label{2.2.1}
\end{eqnarray}
which satisfies with the conservative law (\ref{1.5}). And also, the corresponding numerical finite difference scheme is simplified as follows,
\begin{eqnarray}
&&\frac{v^{n+1}-v^n}{\Delta t}-\frac{\alpha}{2}\Delta_h (u^{n+1}+u^n)=\phi \frac{\cos( u^{n+1})-\cos (u^n)}{u^{n+1}-u^n},\label{2.2.2}\\
&&\frac{u^{n+1}-u^{n}}{\Delta t}=\frac{v^{n+1}+v^n}{2}.   \label{2.2.3}
\end{eqnarray}

Next, we introduce the $l^2$-norm and the  $l^2$ inner product. For any two homogeneous (or periodic) grid functions $f$ and $g$, define the discrete $l^2$ inner product and the discrete $l^2$-norm, respectively, as
\begin{eqnarray}
\langle f,g \rangle =h^2\sum_{i,j=0}^M f_{i,j}g_{i,j},\quad \|f\|_2=\sqrt{\langle f,f \rangle }\label{2.2.4},
\end{eqnarray}
and the following summation by parts is also straightforward,
\begin{eqnarray}
\langle \Delta_h f,  g \rangle =-\langle \nabla_h f,\nabla_h g\rangle, \label{2.2.5}
\end{eqnarray}
with
\begin{eqnarray}
&&\|\nabla_h f\|^2_2=\|D_x f\|^2_2+\|D_y f\|^2_2, \nonumber\\
&&\|D_x f\|^2_2=h^2\sum_{i,j=0}^M (f_{i+1,j}-f_{i,j})^2/h^2,\quad \|D_y f\|^2_2
=h^2\sum_{i,j=0}^M (f_{i,j+1}-f_{i,j})^2/h^2. \label{2.2.6}
\end{eqnarray}

\begin{theorem}\label{thm1}
The scheme (\ref{2.2.2})-(\ref{2.2.3}) is conservative for the  discrete energy, namely,
\begin{eqnarray}
{\cal E}^n=\frac{1}{2}\|v^n\|^2_2+\frac{\alpha}{2}\|\nabla_h u^n\|^2_2+h^2\sum_{i,j=0}^M (\phi (1-\cos(u^n)))_{ij}={\cal E}^0,\label{2.2.7}
\end{eqnarray}
for $n=1,2,\cdots,N$.
\end{theorem}

\begin{remark}\label{rm1}
The main purpose of  taking this discrete energy form is to be in accordance with its continuous
definition (\ref{1.3}). Alternatively, if we delete the constant  in the term $\phi(1-\cos(u^n))$,
the new discrete energy is also conservative.
\end{remark}

\begin{remark}\label{rm2}
In \cite{Asg2013}, the authors had simulated the discrete energy which looks to  likely be conservative. However, it might be caused by the high accuracy of the proposed pseudospectral method. As we know, the explicit numerical scheme could not ensure the conversation generally.
\end{remark}

{\bf Proof.} Taking the inner product of (\ref{2.2.2}) with $u^{n+1}-u^n$ yields
\begin{eqnarray}
&&\langle \frac{v^{n+1}-v^n}{\Delta t},u^{n+1}-u^n \rangle -\frac{\alpha}{2} \langle \Delta_h (u^{n+1}+u^n),u^{n+1}-u^n \rangle
- \langle \phi \frac{\cos( u^{n+1})-\cos (u^n)}{u^{n+1}-u^n}, u^{n+1}-u^n\rangle  \nonumber\\
&&\qquad\qquad\qquad\qquad\qquad  =0.
  \label{2.2.8}
\end{eqnarray}
For the first term, we have
\begin{eqnarray}
\langle \frac{v^{n+1}-v^n}{\Delta t},u^{n+1}-u^n \rangle &=&\langle {v^{n+1}-v^n},\frac{u^{n+1}-u^n}{\Delta t} \rangle =\frac{1}{2}\langle {v^{n+1}-v^n},{v^{n+1}+v^n} \rangle\nonumber\\
&=&\frac{1}{2}(\|v^{n+1}\|_2^2-\|v^{n}\|_2^2),  \label{2.2.9}
\end{eqnarray}
where (\ref{2.2.3}) is applied in the second step. According to (\ref{2.2.2}),  the second term can be analyzed as
\begin{eqnarray}
 -\frac{\alpha}{2} \langle \Delta_h (u^{n+1}+u^n),u^{n+1}-u^n \rangle=\frac{\alpha}{2}(\|\nabla_h u^{n+1}\|^2_2-\|\nabla_h u^{n}\|^2_2).  \label{2.2.10}
\end{eqnarray}
Moreover, for the nonlinear term, we obtain the following result,
\begin{eqnarray}
- \langle \phi \frac{\cos( u^{n+1})-\cos (u^n)}{u^{n+1}-u^n}, u^{n+1}-u^n\rangle =h^2\sum_{i,j=0}^n \phi_{i,j}((1-\cos u^{n+1})-(1-\cos u^n))_{i,j}.  \label{2.2.11}
\end{eqnarray}
By the definition of ${\cal E}^n$, (\ref{2.2.7}) is obtained from (\ref{2.2.8})-(\ref{2.2.11}).

\section{Linear iteration algorithm}
\label{S3}

Since the nonlinear term $\cos(u^{n+1})$ is treated implicitly in (\ref{2.1.8}), it leads to a sequence of nonlinear coupled equations. In order to solve the nonlinear system (\ref{2.1.8}) arising from the implicit treatment, we proposed the following linear iteration algorithm:
\begin{eqnarray}
&&\frac{2u^{n+1,(m+1)}}{\Delta t^2}+\frac{\beta}{\Delta t}u^{n+1,(m+1)}-\frac{\alpha}{2}\Delta_h u^{n+1,(m+1)}\nonumber\\
&&\quad =\phi \frac{\cos( u^{n+1,(m)})-\cos (u^n)}{u^{n+1,(m)}-u^n}+\frac{\beta}{\Delta t}u^{n}+\frac{\alpha}{2}\Delta_h u^{n}+\kappa (u^n,v^n)+F^{n+\frac{1}{2}},\label{2.3.1}
\end{eqnarray}
where $u^{n+1,(m)}$ denotes the approximation solution at the $m$-th iteration.

\begin{theorem}\label{thm2}
The linear iteration scheme (\ref{2.3.1}) is a contraction mapping, provided that
$\Delta t < \frac{\beta+\sqrt{\beta^2+8\phi_0}}{2\phi_0}$.
\end{theorem}

{\bf Proof.} Define the iteration error of each stage via
\begin{eqnarray}
e^{(m)}=u^{n+1,(m)}-u^{n+1}, \label{2.3.2}
\end{eqnarray}
where $u^{n+1,(m)}$ is the $m$-th iteration result generated by the linear iteration
scheme (\ref{2.3.1}). Subtracting (\ref{2.3.1}) from (\ref{2.1.8}) leads to
\begin{eqnarray}
(\frac{2}{\Delta t^2}+\frac{\beta}{\Delta t}-\frac{\alpha}{2}\Delta_h)e^{(m+1)}=\phi \big(\frac{\cos( u^{n+1,(m)})-\cos (u^n)}{u^{n+1,(m)}-u^n}-\frac{\cos( u^{n+1})-\cos (u^n)}{u^{n+1}-u^n}\big ). \label{2.3.3}
\end{eqnarray}
Taking the inner product of (\ref{2.3.3}) with $e^{(m+1)}$, we have
\begin{eqnarray}
&&\langle (\frac{2}{\Delta t^2}+\frac{\beta}{\Delta t}-\frac{\alpha}{2}\Delta_h)e^{(m+1)},e^{(m+1)}\rangle  =(\frac{2}{\Delta t^2}+\frac{\beta}{\Delta t})\|e^{(m+1)}\|^2_2+\frac{\alpha}{2}\|\nabla_h e^{(m+1)}\|^2_2 \nonumber\\
&& =\langle \phi \big(\frac{\cos( u^{n+1,(m)})-\cos (u^n)}{u^{n+1,(m)}-u^n}-\frac{\cos( u^{n+1})
-\cos (u^n)}{u^{n+1}-u^n}\big ),e^{(m+1)}\rangle.  \label{2.3.4}
\end{eqnarray}

Now, we analyze the right-hand side of (\ref{2.3.4}) in detail. For convenience, let
\begin{eqnarray*}
h(x)=\frac{\cos x-\cos a}{x-a}.
\end{eqnarray*}
Using the Lagrange theorem, we obtain $h(x)=-\sin \xi$, where $\xi$ is between $x$ and $a$.
We also compute the derivative of $h(x)$,
\begin{eqnarray}
h'(x)=\frac{-(x-a)\sin x-(\cos x-\cos a)}{(x-a)^2}=\frac{-\sin x-\frac{(\cos x-\cos a)}{x-a}}{x-a}=\frac{-\sin x+\sin \xi}{x-a}.  \label{2.3.5}
\end{eqnarray}
Applying the Lagrange theorem again for $(-\sin x+\sin \xi)$ yields
\begin{eqnarray}
|h'(x)|=|\cos \xi_1|\cdot \frac{|x-\xi|}{|x-a|}<1,  \label{2.3.6}
\end{eqnarray}
where $\xi_1$ is between $x$ and $\xi$, and the fact that $|x-\xi|< |x-a|$ is used.

Going back to (\ref{2.3.4}) and setting $a=u^n$, we have
\begin{eqnarray}
&&|\frac{\cos( u^{n+1,(m)})-\cos (u^n)}{u^{n+1,(m)}-u^n}-\frac{\cos( u^{n+1})-\cos (u^n)}{u^{n+1}-u^n}|\nonumber\\
&&\qquad =|h(u^{n+1,(m)})-h(u^{n+1})|\nonumber\\
&&\qquad =|h'(\xi_2)|\cdot |u^{n+1,(m)}-u^{n+1}|<|e^{(m)}|, \label{2.3.7}
\end{eqnarray}
with $\xi_2$ between $u^{n+1,(m)}$ and $u^{n+1}$. In turn, one can get
\begin{eqnarray}
&&\langle \phi \big(\frac{\cos( u^{n+1,(m)})-\cos (u^n)}{u^{n+1,(m)}-u^n}-\frac{\cos( u^{n+1})-\cos (u^n)}{u^{n+1}-u^n}\big ),e^{(m+1)}\rangle \nonumber\\
&&\qquad \le \phi_0 |\langle e^{(m)}, e^{(m+1)} \rangle| \le \frac{\phi_0}{2}(\|e^{(m+1)}\|^2_2+\|e^{(m)}\|^2_2).  \label{2.3.8}
\end{eqnarray}
in which $\phi_0$ is the upper bound of $\phi$.

As a result, it follows from the combination of (\ref{2.3.4}) and (\ref{2.3.8}) that
\begin{eqnarray}
(\frac{2}{\Delta t^2}+\frac{\beta}{\Delta t}-\frac{\phi_0}{2})\|e^{(m+1)}\|^2_2+\frac{1}{2}\|\nabla_h e^{(m+1)}\|^2_2 \le \frac{\phi_0}{2}\|e^{(m)}\|^2_2.  \label{2.3.9}
\end{eqnarray}
Therefore, the contraction mapping property can be assured if
\begin{eqnarray}
(\frac{2}{\Delta t^2}+\frac{\beta}{\Delta t}-\frac{\phi_0}{2}) > \frac{\phi_0}{2},  \label{2.3.10}
\end{eqnarray}
which shows that the result is proven.

\begin{remark}\label{rm3}
For a high dimensional problem, either the ADI scheme  or the predictor-corrector scheme is often used to implement the implicit finite difference scheme. Both schemes belong to two-step method or multi-step one for reducing the dimension complexity. For 2D sine-Gordon equation, see \cite{Bratsos2006, Bratsos2007, Cui2010}. Without any decompositions for 2D sine-Gordon equation, the implicit scheme can be solved efficiently by the linear iteration. Certainly, this iteration method is also applied in many numerical method, but the contraction condition for iterations is seldom investigated. On the other hand, we also can present the high-order finite difference scheme if more complicated operators are involved, and can solve it by the iteration method.
\end{remark}

\section{ Truncation errors and the convergence analysis}
\label{S4}

\subsection{Truncation errors}

Let $u_e$ and $v_e=\partial_t u_e$ be  exact solutions of the problem (\ref{1.1})-(\ref{1.3}),
then truncation errors of the
scheme (\ref{2.1.1}) are obtained at discrete grid points as follows,
\begin{eqnarray}
&&\frac{v_e^{k+1}-v_e^k}{\Delta t}+\frac{\beta}{\Delta t}(u_e^{k+1}-u_e^k)-\frac{\alpha}{2}\Delta_h (u_e^{k+1}+u_e^k)-\phi \frac{\cos( u_e^{k+1})-\cos (u_e^k)}{u_e^{k+1}-u_e^k}-F^{k+\frac{1}{2}}=\rho^k,\label{3.1.1}\\
&&\frac{u_e^{k+1}-u_e^{k}}{\Delta t}-\frac{v_e^{k+1}+v_e^k}{2}=s^k. \label{3.1.2}
\end{eqnarray}

In fact, for a given function $f(x) \in C^5$, one can get
\begin{eqnarray}
\frac{f(\xi)-f(\eta)}{\xi -\eta}=f'(\frac{\xi +\eta}{2})+\frac{1}{24}f''(\frac{\xi +\eta}{2})(\xi-\eta)^2+O((\xi-\eta)^4).\label{3.1.3}
\end{eqnarray}
Taking $f(x)=\cos(x), \xi=u_e^{k+1}$ and $\eta=u_e^{k}$, we have
\begin{eqnarray}
 \frac{\cos( u_e^{k+1})-\cos (u_e^k)}{u_e^{k+1}-u_e^k}=-\sin(\frac{u_e^{k+1}+u_e^{k}}{2})+O(\Delta t^2)
 =-\sin(u_e^{k+\frac{1}{2}})+O(\Delta t^2).\label{3.1.4}
\end{eqnarray}
Similarly, other terms can be analyzed and the details are omitted.
Hence, we obtain the following lemma.

\begin{lemma}
\label{lm1} Suppose $u_e$ and $v_e$ are smooth enough, then $|\rho^k_{ij}|+|s^k_{ij}|=O(\Delta t^2+h^2)$ holds as $\Delta t, h\rightarrow 0$.
\end{lemma}

\subsection{ Convergence analysis}

Define  discrete error functions as follows,
\begin{eqnarray}
\ti u^k=u_e^k-u^k, \qquad \ti v^k=v_e^k-v^k.\label{3.2.1}
\end{eqnarray}

Now, we present the following convergence result.
\begin{theorem}\label{thm3}
Assume that $u_e$ and $v_e=\partial_t u_e$ are the exact solutions of the problem (\ref{1.1})-(\ref{1.3}),
and denote $(u,v)$ as the numerical solution given by the finite difference scheme (\ref{2.1.3}) and (\ref{2.1.4}). Then, we have
\begin{eqnarray}
\|\ti v^n\|^2_2+\|\ti u^n\|_2^2+\alpha \|\nabla_h \ti u^n\|_2^2 \le \ti C \cdot O(\Delta t^2+ h^2),\label{3.2.2}
\end{eqnarray}
where the constant $\ti C$ is dependent on the final time $T$ and is independent
on  $\Delta t$ and $h$.
\end{theorem}

{\bf Proof.} Subtracting (\ref{2.1.1})-(\ref{2.1.2}) from (\ref{3.1.1})-(\ref{3.1.2}),
respectively, we get
\begin{eqnarray}
&&\frac{\ti v^{k+1}-\ti v^k}{\Delta t}+\frac{\beta}{\Delta t}(\ti u^{k+1}-\ti u^k)-\frac{\alpha}{2}\Delta_h (\ti u^{k+1}+\ti u^k)-\phi \big (\frac{\cos( u_e^{k+1})-\cos (u_e^k)}{u_e^{k+1}-u_e^k}-\frac{\cos( u^{k+1})-\cos (u^k)}{u^{k+1}-u^k}\big )\nonumber\\
&&\qquad\qquad\qquad =\rho^k,\label{3.2.3}\\
&&\frac{\ti u^{k+1}-\ti u^{k}}{\Delta t}-\frac{\ti v^{k+1}+\ti v^k}{2}=s^k. \label{3.2.4}
\end{eqnarray}
Taking the inner product of (\ref{3.2.3}) with $\ti u^{k+1}-\ti u^{k}$ yields
\begin{eqnarray}
&&\frac{1}{\Delta t}\langle \ti v^{k+1}-\ti v^k,\ti u^{k+1}-\ti u^{k} \rangle+\frac{\beta}{\Delta t}\langle \ti u^{k+1}-\ti u^k,\ti u^{k+1}-\ti u^k \rangle -\frac{\alpha}{2}\langle \Delta_h (\ti u^{k+1}+\ti u^k), \ti u^{k+1}-\ti u^{k}\rangle \nonumber\\
 &&=\langle \phi \big (\frac{\cos( u_e^{k+1})-\cos (u_e^k)}{u_e^{k+1}-u_e^k}-\frac{\cos( u^{k+1})-\cos (u^k)}{u^{k+1}-u^k}\big ),\ti u^{k+1}-\ti u^{k}\rangle +\langle \rho^k, \ti u^{k+1}-\ti u^{k} \rangle. \label{3.2.5}
\end{eqnarray}

Next, we begin to analyze the nonlinear term on the right-hand of (\ref{3.2.5}).
Noting that
\begin{eqnarray}
 &&|\big (\frac{\cos( u_e^{k+1})-\cos (u_e^k)}{u_e^{k+1}-u_e^k}-\frac{\cos( u^{k+1})-\cos (u^k)}{u^{k+1}-u^k}\big )|  \nonumber\\
 &&\le |\big (\frac{\cos( u_e^{k+1})-\cos (u_e^k)}{u_e^{k+1}-u_e^k}-\frac{\cos( u^{k+1})-\cos (u_e^k)}{u^{k+1}-u_e^k}\big )| \nonumber\\
 &&\qquad + |\big (\frac{\cos( u^{k+1})-\cos (u_e^k)}{u^{k+1}-u_e^k}-\frac{\cos( u^{k+1})-\cos (u^k)}{u^{k+1}-u^k}\big )| \label{3.2.6}
\end{eqnarray}
and recalling the definition (\ref{2.3.5}) of $h(x)$, we have the following results
\begin{eqnarray}
&&|\big (\frac{\cos( u_e^{k+1})-\cos (u_e^k)}{u_e^{k+1}-u_e^k}-\frac{\cos( u^{k+1})-\cos (u_e^k)}{u^{k+1}-u_e^k}\big )|\nonumber\\
&&\qquad \le |h(u_e^{k+1})-h(u^{k+1})|=|h'(\eta_1)|\cdot |u_e^{k+1}-u^{k+1}|\le |\ti u^{k+1}|,   \label{3.2.7}
\end{eqnarray}
in which we choose $a=u_e^k$, and
\begin{eqnarray}
&&|\big (\frac{\cos( u^{k+1})-\cos (u_e^k)}{u^{k+1}-u_e^k}-\frac{\cos( u^{k+1})-\cos (u^k)}{u^{k+1}-u^k}\big )|\nonumber\\
&&\qquad \le |h(u_e^{k})-h(u^{k})|=|h'(\eta_2)|\cdot |u_e^{k}-u^{k}|\le |\ti u^{k}|,   \label{3.2.8}
\end{eqnarray}
in which  $a=u^{k+1}$, respectively. Substituting (\ref{3.2.7}) and (\ref{3.2.8}) into (\ref{3.2.6}) gives
\begin{eqnarray}
 |\big (\frac{\cos( u_e^{k+1})-\cos (u_e^k)}{u_e^{k+1}-u_e^k}-\frac{\cos( u^{k+1})-\cos (u^k)}{u^{k+1}-u^k}\big )| \le |\ti u^{k+1}|+|\ti u^{k}|.
 \label{3.2.9}
\end{eqnarray}

Also, it follows from (\ref{3.2.4}) that
\begin{eqnarray*}
\ti u^{k+1}-\ti u^{k}=\frac{\Delta t}{2}(\ti v^{k+1}+\ti v^k)+\Delta t s^k.
\end{eqnarray*}
Therefore, we can obtain the estimate as follows,
\begin{eqnarray}
&&\langle \phi \big (\frac{\cos( u_e^{k+1})-\cos (u_e^k)}{u_e^{k+1}-u_e^k}-\frac{\cos( u^{k+1})-\cos (u^k)}{u^{k+1}-u^k}\big ),\ti u^{k+1}-\ti u^{k}\rangle \nonumber\\
&&\le \phi_0 \langle |\ti u^{k+1}|+ |\ti u^{k}|, |\ti u^{k+1}-\ti u^{k}|\rangle \nonumber\\
&&\le \frac{\phi_0}{2} \Delta t \langle |\ti u^{k+1}|+ |\ti u^{k}|, |\ti v^{k+1}|+|\ti v^{k}|\rangle +\phi_0 \Delta t \langle |\ti u^{k+1}|+ |\ti u^{k}|, |s^k|\rangle \nonumber\\
&&\le \frac{\phi_0}{4} \Delta t \big ( \|(|\ti u^{k+1}|+ |\ti u^{k}|)\|_2^2+\|(|\ti v^{k+1}|+ |\ti v^{k}|)\|_2^2 \big )+\frac{\phi_0}{2} \Delta t \big ( \|(|\ti u^{k+1}|_2^2+ |\ti u^{k}|)\|_2^2 \big ) +\phi_0 \Delta t \|s^k\|_2^2 \nonumber\\
&&\le \frac{\phi_0}{2} \Delta t \big ( \|\ti u^{k+1}\|_2^2+ \|\ti u^{k}\|_2^2+\|\ti v^{k+1}\|_2^2+ \|\ti v^{k}\|_2^2 \big )+\phi_0 \Delta t \big ( \|\ti u^{k+1}\|_2^2+ \|\ti u^{k}\|_2^2 \big ) +\phi_0 \Delta t \|s^k\|_2^2 \nonumber\\
&&\le \frac{3\phi_0}{2} \Delta t \big ( \|\ti u^{k+1}\|_2^2+ \|\ti u^{k}\|_2^2\big )+\frac{\phi_0}{2} \Delta t \big (\|\ti v^{k+1}\|_2^2+ \|\ti v^{k}\|_2^2 \big ) +\phi_0 \Delta t \|s^k\|_2^2,
\label{3.2.10}
\end{eqnarray}
where the inequality $(a+b)^2 \le 2(a^2+b^2)$ is used.\\

For the first term in (\ref{3.2.5}), using (\ref{3.2.4}), we can arrive at,
\begin{eqnarray}
&&\frac{1}{\Delta t}\langle \ti v^{k+1}-\ti v^k,\ti u^{k+1}-\ti u^{k} \rangle
\nonumber\\
&&\qquad =\frac{1}{2}\langle \ti v^{k+1}-\ti v^k,\ti v^{k+1}+\ti v^k \rangle
+\Delta t \langle \ti v^{k+1}-\ti v^k,s^k \rangle  \nonumber\\
&&\qquad  =\frac{1}{2} \big (  \|\ti v^{k+1}\|_2^2- \|\ti v^{k}\|_2^2  \big )+ \Delta t \langle \ti v^{k+1}-\ti v^k,s^k \rangle   \nonumber\\
&&\qquad  \ge  \frac{1}{2} \big (  \|\ti v^{k+1}\|_2^2- \|\ti v^{k}\|_2^2  \big )- \frac{1}{2} \Delta t \big (  \|\ti v^{k+1}-\ti v^{k} \|_2^2+\|s^k\|_2^2 \big ) \nonumber\\
&&\qquad   \ge \frac{1}{2} \big (  \|\ti v^{k+1}\|_2^2- \|\ti v^{k}\|_2^2  \big )-  \Delta t \big (  \|\ti v^{k+1}\|^2_2+\|\ti v^{k} \|_2^2\big )- \frac{1}{2} \Delta t \|s^k\|_2^2.
\label{3.2.11}
\end{eqnarray}
And it follows from the second term in (\ref{3.2.5}) that
\begin{eqnarray}
\frac{\beta}{\Delta t}\langle \ti u^{k+1}-\ti u^k,\ti u^{k+1}-\ti u^k \rangle=\frac{\beta}{\Delta t}\|\ti u^{k+1}-\ti u^k \|_2^2.  \label{3.2.12}
\end{eqnarray}

As for the diffusion term,  the following result is also straightforward from (\ref{2.2.2}),
\begin{eqnarray}
-\langle \Delta_h (\ti u^{k+1}+\ti u^{k}),  \ti u^{k+1}-\ti u^{k} \rangle &=&\langle \nabla_h (\ti u^{k+1}+\ti u^{k}),\nabla_h (\ti u^{k+1}-\ti u^{k})\rangle, \nonumber\\
&=&\|\nabla_h \ti u^{k+1}\|_2^2-\|\nabla_h \ti u^{k}\|_2^2. \label{3.2.13}
\end{eqnarray}
In turn, considering the local truncation error term in (\ref{3.2.5}), we have
\begin{eqnarray}
\langle \rho^k,  \ti u^{k+1}-\ti u^{k} \rangle &=& \Delta t \langle \rho^k,  \frac{\ti u^{k+1}-\ti u^{k}}{\Delta t} \rangle =  \frac{1}{2} \Delta t \langle \rho^k,  \ti v^{k+1}+\ti v^{k} \rangle + \Delta t   \langle \rho^k,  s^k \rangle \nonumber\\
&\le& \frac{1}{4}\Delta t (\|\rho^k\|_2^2+\|\ti v^{k+1}+\ti v^{k}\|_2^2)+\frac{1}{2}\Delta t (\|\rho^k\|_2^2+\|s^k\|_2^2) \nonumber\\
&\le& \frac{3}{4}\Delta t \|\rho^k\|_2^2 +\frac{1}{2}\Delta t \|s^k\|_2^2+\frac{1}{2}\Delta t (\|\ti v^{k+1}\|_2^2+\|\ti v^{k}\|_2^2).
 \label{3.2.14}
\end{eqnarray}

Consequently, it follows from (\ref{3.2.5}), (\ref{3.2.11}), (\ref{3.2.12}), (\ref{3.2.13}) and (\ref{3.2.14}) that
\begin{eqnarray}
&&\frac{1}{2} \big (  \|\ti v^{k+1}\|_2^2- \|\ti v^{k}\|_2^2\big )  +\frac{\alpha}{2}\big ( \|\nabla_h \ti u^{k+1}\|_2^2-\|\nabla_h \ti u^{k}\|_2^2 \big )+\frac{\beta}{\Delta t}\|\ti u^{k+1}-\ti u^k \|_2^2\nonumber\\
&&\le \frac{3}{2}\phi_0 \Delta t(\|\ti u^{k+1}\|_2^2+\|\ti u^{k}\|_2^2)+(\frac{3}{2}+\frac{\phi_0}{2})\Delta t(\|\ti v^{k+1}\|_2^2+\|\ti v^{k}\|_2^2)\nonumber\\
&&\qquad\qquad\qquad\qquad +\frac{3}{4} \Delta t \|\rho^k\|_2^2+(\phi_0 +1)\Delta t\|s^k\|_2^2.
 \label{3.2.15}
\end{eqnarray}

Unfortunately, we so far can not obtain the estimate for $\|\ti u^k\|$ and $\|\ti v^k\|$ by the discrete Gronwall inequality, since the term $\|\ti u^k\|$ does not appear in the left-hand part of (\ref{3.2.15}). To remedy this, we need to analyze the estimate for $\|\ti u^k\|$.

Going back to (\ref{3.2.4}), it  can be rewritten as
\begin{eqnarray}
\ti u^{k}-\ti u^{k-1}=\frac{1}{2} \Delta t (\ti v^{k}+\ti v^{k-1})+ \Delta t s^{k-1}. \label{3.2.16}
\end{eqnarray}
Summing up (\ref{3.2.16}) from 1 to $k$, we get
\begin{eqnarray}
\ti u^{k}=\ti u^{0}+\frac{1}{2} \Delta t \sum_{l=1}^k(\ti v^{l}+\ti v^{l-1})+ \Delta t \sum_{l=1}^{k-1} s^l =\frac{1}{2} \Delta t \sum_{l=1}^k(\ti v^{l}+\ti v^{l-1})+ \Delta t \sum_{l=1}^{k-1} s^l, \label{3.2.17}
\end{eqnarray}
where we apply the fact that $\ti u^{0}=0$. In turn, an application of Cauchy inequality implies that
\begin{eqnarray}
|\ti u^{k}|^2 &\le& 2\big [\frac{1}{4} \Delta t^2 (\sum_{l=1}^k(\ti v^{l}+\ti v^{l-1}))^2+ \Delta t^2 (\sum_{l=1}^{k-1} s^l)^2 \big ]\nonumber\\
&\le& \frac{1}{2} \Delta t^2 \sum_{l=1}^k 2(|\ti v^{l}|^2+|\ti v^{l-1}|^2)+ \Delta t^2(\sum_{l=1}^{k-1} O(\Delta t^2))^2  \nonumber\\
&\le& 2k \Delta t^2  \sum_{l=1}^k |\ti v^{l}|^2 + k^2 \Delta t^2( O(\Delta t^2))^2 \nonumber\\
&\le& 2T \Delta t \sum_{l=1}^k |\ti v^{l}|^2 + T^2 O(\Delta t^4), \label{3.2.18}
\end{eqnarray}
in which $k\Delta t\le T $ is used in the last step. This shows
\begin{eqnarray}
|\ti u^{k}|^2 \le 2T \Delta t \sum_{l=1}^k |\ti v^{l}|^2 + C\cdot O(\Delta t^4), \quad |\ti u^{k+1}|^2 \le 2T \Delta t \sum_{l=1}^{k+1} |\ti v^{l}|^2 + C \cdot O(\Delta t^4). \label{3.2.19}
\end{eqnarray}

Substituting the result above into (\ref{3.2.15}) leads to
\begin{eqnarray}
&&\frac{1}{2} \big (  \|\ti v^{k+1}\|_2^2- \|\ti v^{k}\|_2^2\big )  +\frac{\alpha}{2}\big ( \|\nabla_h \ti u^{k+1}\|_2^2-\|\nabla_h \ti u^{k}\|_2^2 \big )+\frac{\beta}{\Delta t}\|\ti u^{k+1}-\ti u^k \|_2^2 \nonumber\\
&&\le 6\phi_0 T \Delta t^2\sum_{l=0}^{k+1}\|\ti v^{l}\|_2^2+(\frac{3}{2}+\frac{\phi_0}{2})\Delta t(\|\ti v^{k+1}\|_2^2
+\|\ti v^{k}\|_2^2)+C\Delta t \cdot O(\Delta t^4+h^4). \label{3.2.20}
\end{eqnarray}
Summing over $k$ from 0 to $n-1$, and using $\ti u^{0}=0, \ti v^{0}=0 $, we have
\begin{eqnarray}
\frac{1}{2}  \|\ti v^{n}\|_2^2  + \frac{\alpha}{2}\|\nabla_h \ti u^{n}\|_2^2
\le 6\phi_0 T \Delta t^2\sum_{k=0}^{n-1}\sum_{l=0}^{k+1}\|\ti v^{l}\|_2^2+(3+\phi_0)\Delta t \sum_{k=0}^{n-1}\|\ti v^{k}\|_2^2
+CT \cdot O(\Delta t^4+h^4). \label{3.2.21}
\end{eqnarray}
Then according to  $\sum_{l=0}^{k+1}\|\ti v^{l}\|_2^2 \le \sum_{l=0}^{n}\|\ti v^{l}\|_2^2$ for $\forall \quad 0\le k+1 \le n$, the inequality above can be rewritten as
\begin{eqnarray}
\frac{1}{2}  \|\ti v^{n}\|_2^2  + \frac{\alpha}{2}\|\nabla_h \ti u^{n}\|_2^2
\le (6\phi_0 T^2+3+\phi_0) \Delta t\sum_{k=0}^{n}\|\ti v^{k}\|_2^2+CT \cdot O(\Delta t^4+h^4). \label{3.2.22}
\end{eqnarray}
Let $\ti E^n=\|\ti v^{n}\|_2^2  +\alpha \|\nabla_h \ti u^{n}\|_2^2 $. Then we  get
\begin{eqnarray}
\ti E^n \le (12\phi_0 T^2+6+2\phi_0) \Delta t\sum_{l=0}^{n-1}\ti E^{l}+CT \cdot O(\Delta t^4+h^4). \label{3.2.23}
\end{eqnarray}

By the discrete Gronwall inequality, we derive that
\begin{eqnarray}
\|\ti v^{n}\|_2^2  +\alpha \|\nabla_h \ti u^{n}\|_2^2  \le \ti C \cdot O(\Delta t^4+h^4), \label{3.2.24}
\end{eqnarray}
that is,
\begin{eqnarray}
\|\ti v^{n}\|_2  + \alpha \|\nabla_h \ti u^{n}\|_2  \le \ti C \cdot O(\Delta t^2+h^2), \label{3.2.25}
\end{eqnarray}
where $\ti C $ is independent on $\Delta t$ and $h$.

Moreover, from (\ref{3.2.19}), we have
\begin{eqnarray}
|\ti u^{n}|^2 \le 2T \Delta t \sum_{l=1}^k |\ti v^{l}|^2 + C\cdot O(\Delta t^4)\le C \cdot O(\Delta t^4+h^4). \label{3.2.26}
\end{eqnarray}
Finally, the combination of (\ref{3.2.25}) and (\ref{3.2.26}) gives
\begin{eqnarray*}
\|\ti v^{n}\|_2 +\|\ti u^{n}\|_2 + \alpha \|\nabla_h \ti u^{n}\|_2  \le \ti C \cdot O(\Delta t^2+h^2),
\end{eqnarray*}
which shows the the unconditional convergence in  the sense of $l^2$-norm is obtained.  This completes the proof of Theorem \ref{thm3}.

\begin{remark}\label{rm4}
For the one dimensional sine-Gordon equation, by the discrete version of the Sobolev imbedding inequality,
the $\|\cdot\|_{\infty}$ estimate can be obtained from the conservative property and the $l^2$-norm. But for
two dimension case, $\|\cdot\|_{\infty}$ error bounds of the numerical scheme (\ref{2.1.1}) and (\ref{2.1.2})
maybe slightly complicated in obtaining the a priori uniform estimate of the numerical solution. It also can
be obtained if we perform a higher consistency analysis by a careful Taylor expansion. The details are
skipped for simplicity of presentation and an analogous technique can be seen in \cite{Wang2015}.
\end{remark}

\section{Numerical simulations}
\label{S5}

\setcounter{equation}{0}

We now perform a couple of numerical experiments that support the theoretical results and error estimates
for the scheme given by (\ref{2.1.3})-(\ref{2.1.4}).

\subsection{Verification of the second order accuracy }

\begin{figure}
\centering
\begin{subfigure}[b]{7cm}
\frame{\includegraphics[width=7cm]{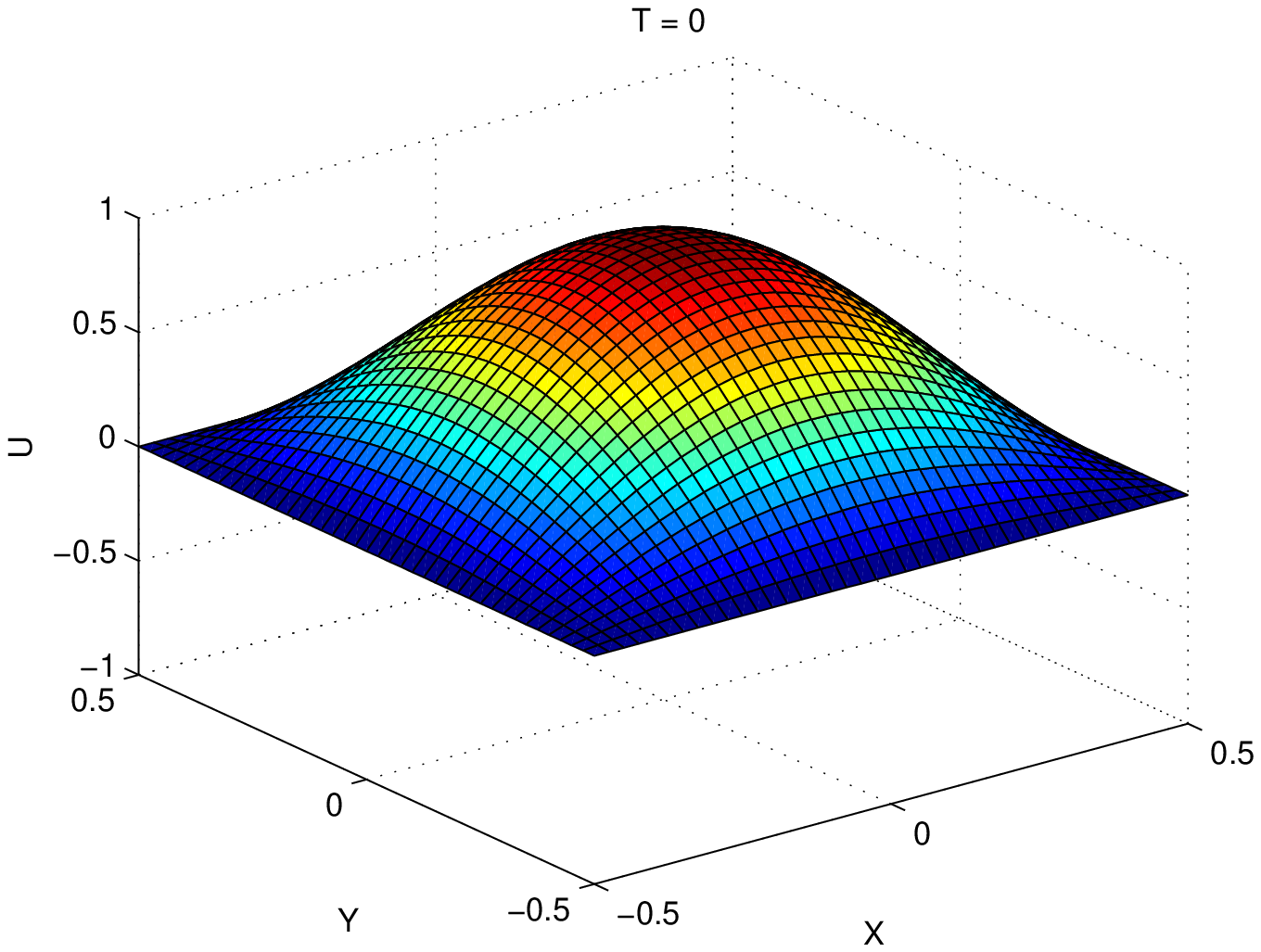}}
\end{subfigure}
\hspace{1cm}
\begin{subfigure}[b]{7cm}
\centering
\frame{\includegraphics[width=7cm]{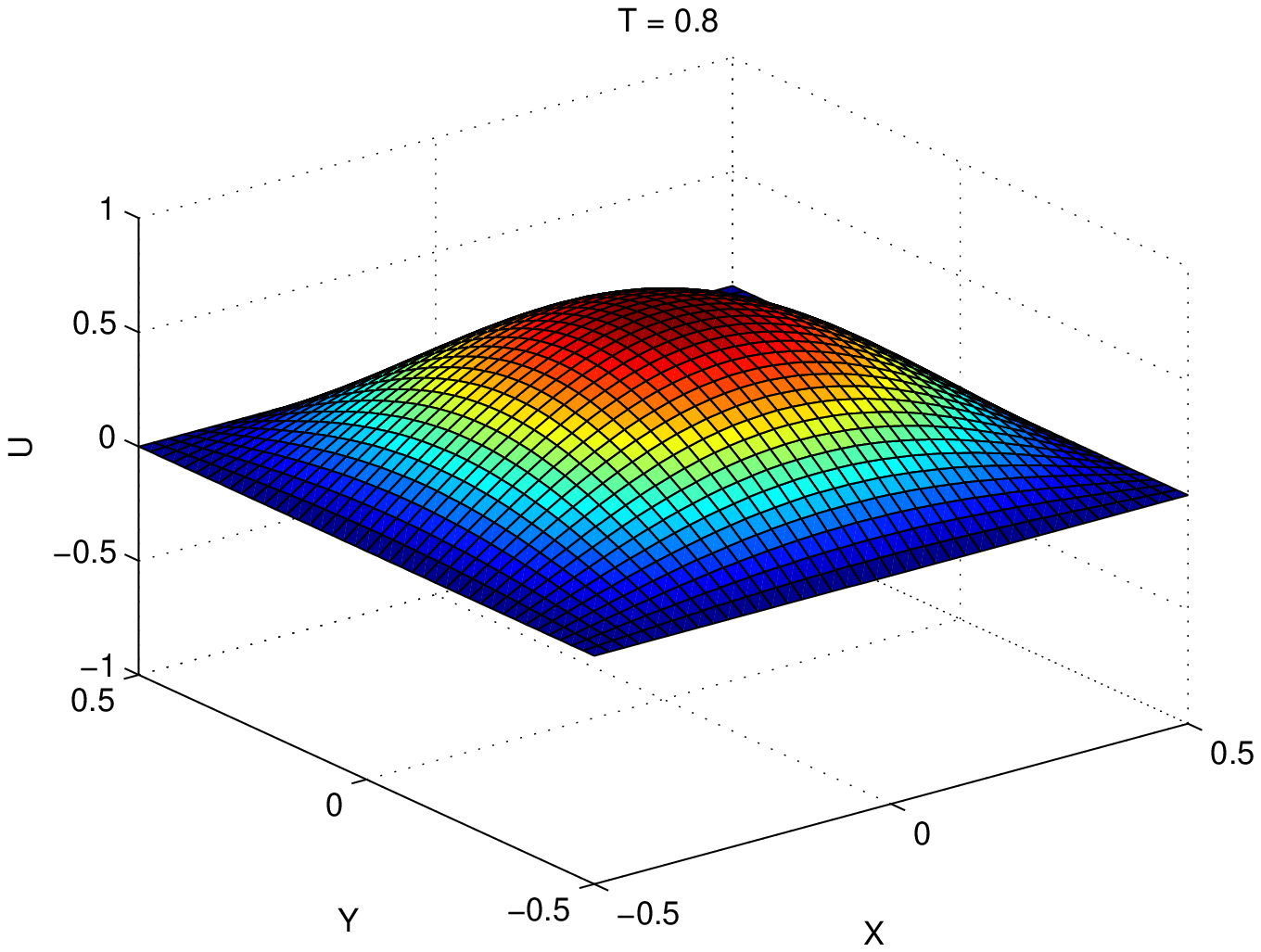}}
\end{subfigure}
\caption{Left: Profile of the exact solution  $u(x, y,t)$ at $t=0$. Right: Numerical solutions of $u$ at $T=0.8$ with $h=0.025$ and $\Delta t=0.1$.}
\label{Fig12}
\end{figure}

\begin{figure}
\centering
\begin{subfigure}[b]{7cm}
\frame{\includegraphics[width=7cm]{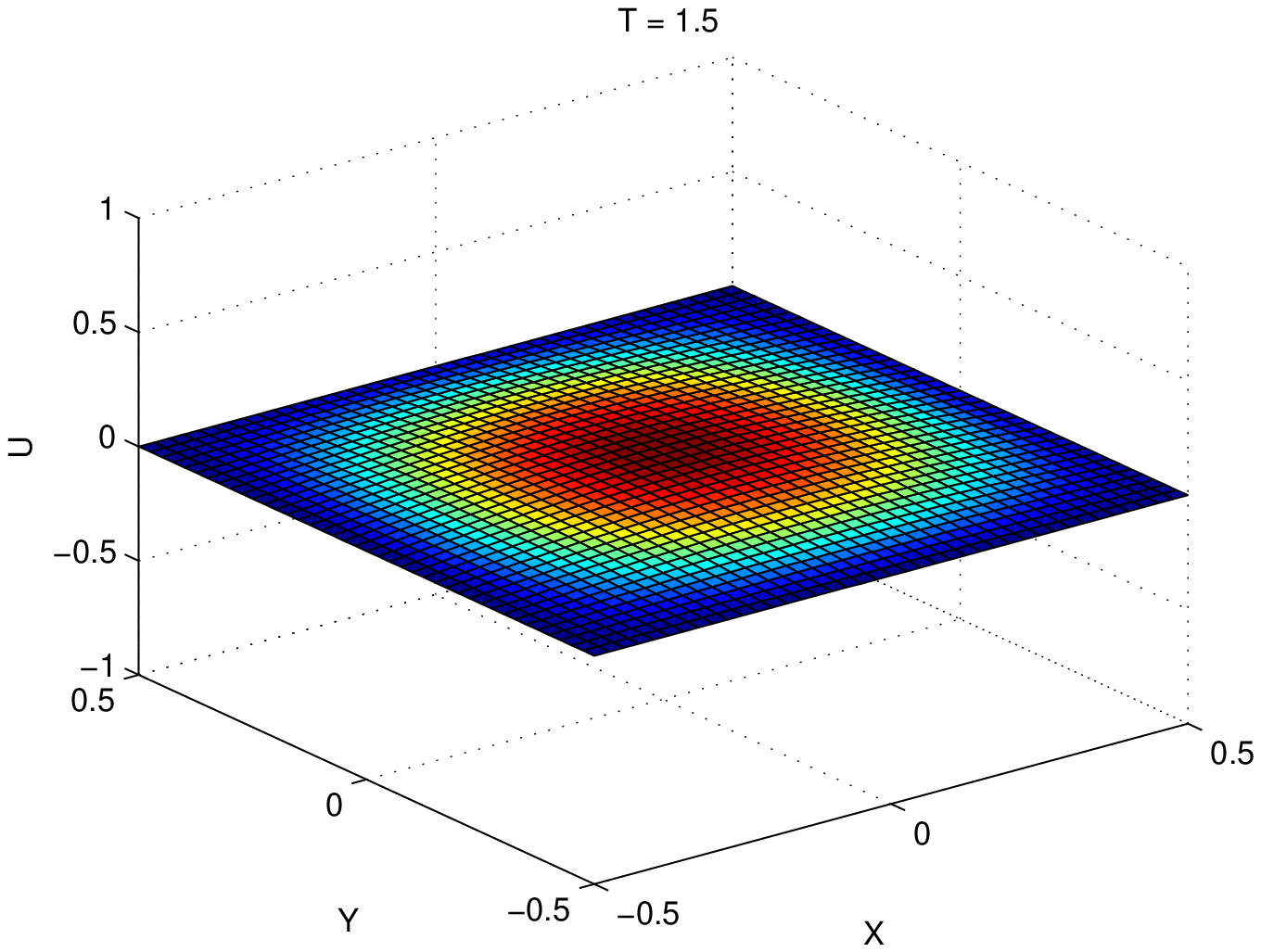}}
\end{subfigure}
\hspace{1cm}
\begin{subfigure}[b]{7cm}
\centering
\frame{\includegraphics[width=7cm]{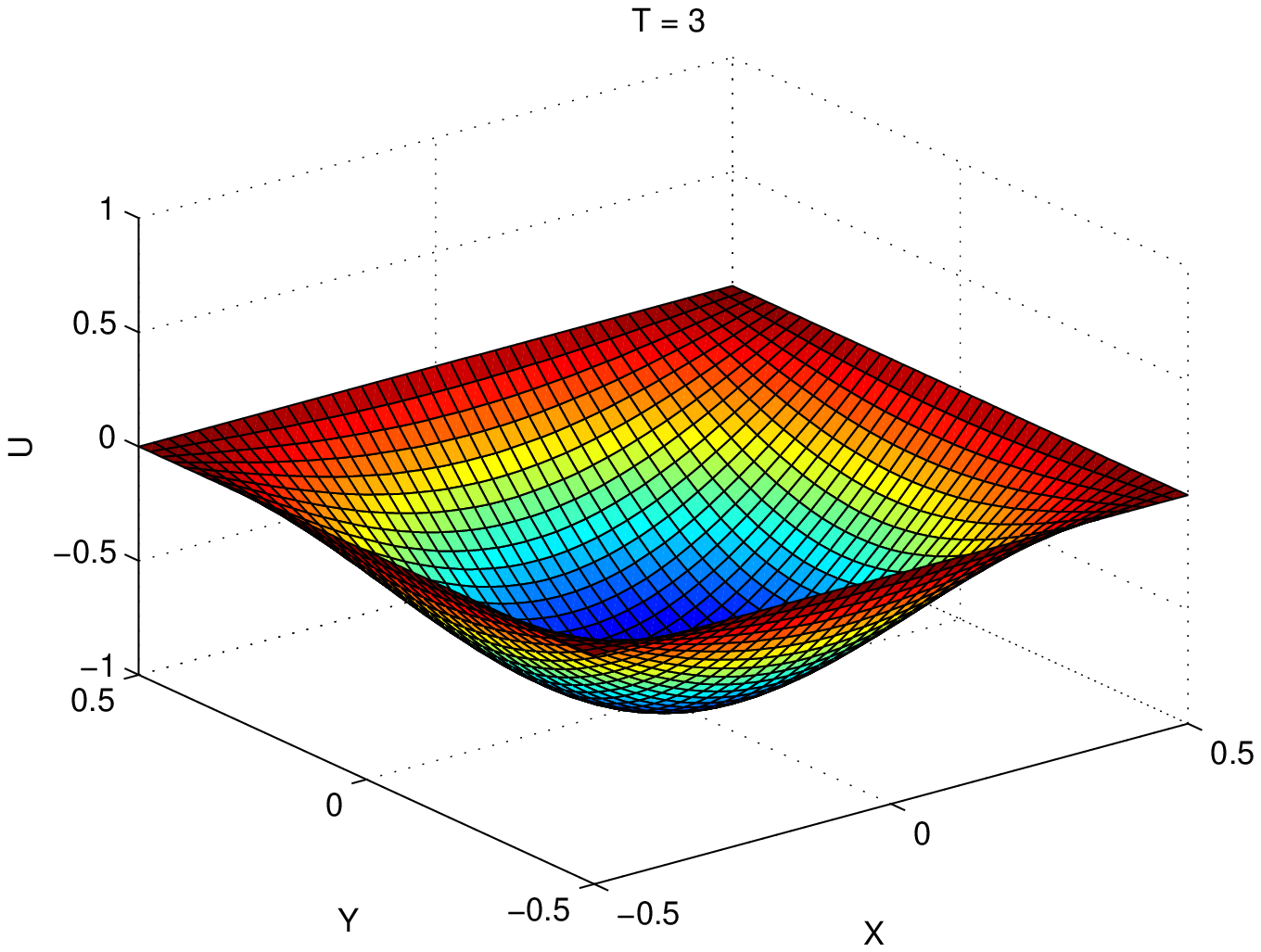}}
\end{subfigure}
\caption{Left: Numerical solutions of $u$ at $T=1.5$ with $h=0.025$ and $\Delta t=0.1$. Right: Numerical solutions of $u$ at $T=3$ with $h=0.025$ and $\Delta t=0.1$.}
\label{Fig34}
\end{figure}

In the first test, we consider the equation (\ref{1.1}) in the  domain $\Omega=[-\frac{1}{2},\frac{1}{2}]\times [-\frac{1}{2},\frac{1}{2}]$ with $\beta=0$, $\phi(x, y) = 1$, $\alpha=\frac{1}{2\pi^2}$ and $f(x,y,t)=\sin(\cos(\pi x)\cos(\pi y)\cos(t))$. The exact solution of (\ref{1.1}) is given by
\begin{eqnarray}
&&u_e(x,y,t)=\cos(\pi x)\cos(\pi y)\cos(t), \label{5.1.1}\\
&&v_e(x,y,t)=-\cos(\pi x)\cos(\pi y)\sin(t).\label{5.1.2}
\end{eqnarray}
The  initial conditions and the boundary condition can be obtained from the exact solution.
 Fig. \ref{Fig12} and Fig. \ref{Fig34} show the
 profile of the exact solution $u$ at $t=0$ on $\Omega$, and  the numerical solutions $u^n$
 at time $T =0.8, 1.5, 3$ with $\Delta t=0.1$ and $h=0.025$, respectively.
The errors in the sense of $l^2$-norm of the numerical solutions for different mesh steps
$h$  and $\Delta t$ at time $t =1, 2, 3, 4$ and $5$ can be found in Table \ref{tab1}, and the corresponding
numerical order of convergence are listed in Table \ref{tab2}. Clearly, it verifies the  second order accuracy in Theorem \ref{thm3}.

\begin{table}
\begin{center}
\caption{ The $\|\cdot\|_2$ errors estimates of  $u^n$ and $v^n$ with various values $h$ and $\Delta t$.}
\begin{tabular}{ccccccc} \hline
\multicolumn{1}{c}{}&\multicolumn{2}{c}{$\Delta t=0.2,h=0.1$}&\multicolumn{2}{c}{$\Delta t=0.1,h=0.05$}&\multicolumn{2}{c}{$\Delta t=0.05,h=0.025$}\\\hline
{}&{$\|\ti u^n\|_2$}&{$\|\ti v^n\|_2$}&{$\|\ti u^n\|_2$}&{$\|\ti v^n\|_2$}&{$\|\ti u^n\|_2$}&{$\|\ti v^n\|_2$}\\\hline
t=1 & 3.6332417e-3   &  4.1077748e-3   & 9.1807718e-4   & 1.0275471e-3    & 2.3013872e-4 &  2.5692078e-4  \\
t=2 & 5.3736869e-3   &  5.7848266e-3   & 1.3423090e-3   & 1.4826793e-3    & 3.3549256e-4 &  3.7298456e-4  \\
t=3 & 4.8006248e-3   &  1.4044778e-2   & 1.2425939e-3   & 3.5219515e-3    & 3.1334248e-4 &  8.8112076e-4  \\
t=4 & 1.5080860e-2   &  1.8066820e-3   & 3.7880320e-3   & 3.6020681e-4    & 9.4806730e-4 &  8.4218108e-5  \\
t=5 & 6.3731025e-3   &  2.1098154e-2   & 1.5032623e-3   & 5.3465660e-3    & 3.7005127e-4 &  1.3410308e-3  \\
\hline
\end{tabular}
\label{tab1}
\end{center}
\end{table}

\begin{table}
\begin{center}
\caption{Numerical verification of theoretical accuracy $O(\Delta t^2+h^2).$ }
\begin{tabular}{ccccccc} \hline
\multicolumn{1}{c}{}&\multicolumn{3}{c}{$\|\ti u^n(h,\Delta t)\|_2\big / \|\ti u^{2n}(\frac{h}{2},\frac{\Delta t}{2})\|_2$}&\multicolumn{3}{c}{$\|\ti v^n(h,\Delta t)\|_2 \big / \|\ti v^{2n}(\frac{h}{2},\frac{\Delta t}{2})\|_2$}\\\hline
{}&{$\Delta t=0.2$}&{$\Delta t=0.1$}&{$\Delta t=0.05$}&{$\Delta t=0.2$}&{$\Delta t=0.1$}&{$\Delta t=0.05$}\\
{}&{$h=0.1$}&{$h=0.05$}&{$h=0.025$}&{$h=0.1$}&{$h=0.05$}&{$h=0.025$}\\\hline
t=1 & - &  1.9846  &  1.9961  &-  &  1.9992 &  1.9998 \\
t=2 & - &  2.0012  &  2.0004  &-  &  1.9641 & 1.9910  \\
t=3 & - &  1.9499  &  1.9875  &-  &  1.9956 & 1.9990  \\
t=4 & - &  1.9932  &  1.9984  &-  &  2.3264 & 2.0966  \\
t=5 & - &  2.0839  &  2.0223  &-  &  1.9804 & 1.9953  \\
\hline
\end{tabular}
\label{tab2}
\end{center}
\end{table}

\subsection{Energy conservation for the undamped equation}

In the second test, we consider the  homogeneous boundary condition $u|_{\partial \Omega}=0$ for the equation (\ref{1.4}) on $\Omega=[0,1]\times[0,1]$ with $\phi=1$,
\begin{eqnarray}
&&\varphi_1(x,y)=\sin(2\pi x)\sin(2\pi y), \label{5.2.1}\\
&&\varphi_2(x,y)=0. \label{5.2.2}
\end{eqnarray}
 We take $T=1$ and  the discrete
energy ${\cal E}^n$ at the different time for $\Delta t=0.002, h=0.05$ and $\Delta t=0.001,h=0.025$ can be  found in
Table \ref{tab3}. Obviously, values of ${\cal E}^n$ at the different  time
 remain nearly a constant as time increases.

\begin{table}
\begin{center}
\caption{ The discrete energy ${\cal E}^n$  with different $\Delta t$ and $h$.}
\begin{tabular}{ccccc} \hline
\multicolumn{1}{c}{}&\multicolumn{1}{c}{$\Delta t=0.002,h=0.05$}&\multicolumn{1}{c}{$\Delta t=0.001,h=0.025$}\\\hline
t=0.2  & 10.351697 & 10.087760     \\
t=0.4  & 10.216081 & 10.031903     \\
t=0.6  & 10.151906 & 10.112902     \\
t=0.8  & 10.288459 & 10.199227     \\
t=1    & 10.419025 & 10.280056     \\
\hline
\end{tabular}
\label{tab3}
\end{center}
\end{table}

\subsection{Circular ring soliton}

\begin{figure}
\centering
\begin{subfigure}[b]{7cm}
\frame{\includegraphics[width=7cm]{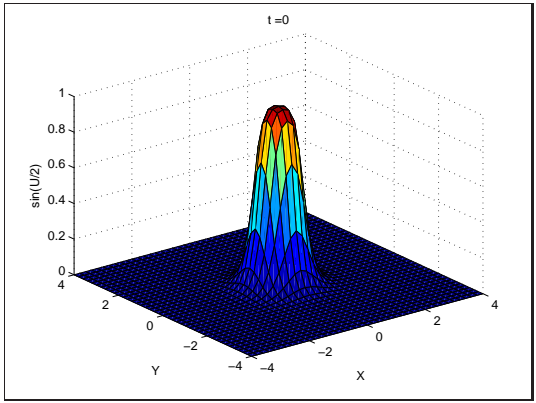}}
\end{subfigure}
\hspace{1cm}
\begin{subfigure}[b]{7cm}
\centering
\frame{\includegraphics[width=7cm]{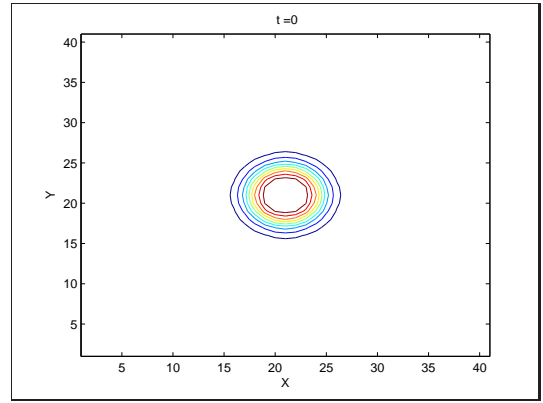}}
\end{subfigure}
\caption{The initial function of $u$ and its contour profile.}
\label{Fig56}
\end{figure}

\begin{figure}
\centering
\begin{subfigure}[b]{7cm}
\frame{\includegraphics[width=7cm]{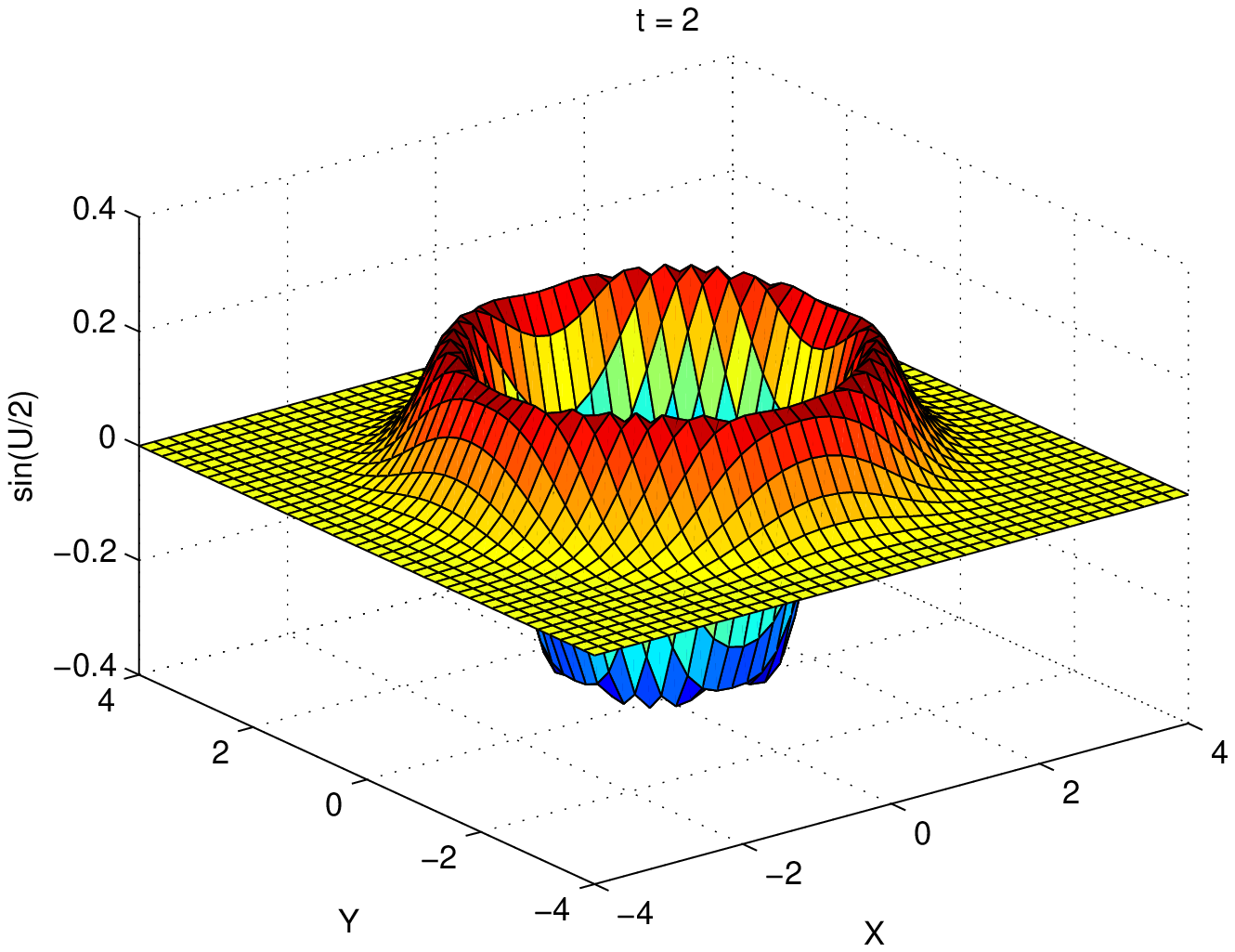}}
\end{subfigure}
\hspace{1cm}
\begin{subfigure}[b]{7cm}
\centering
\frame{\includegraphics[width=7cm]{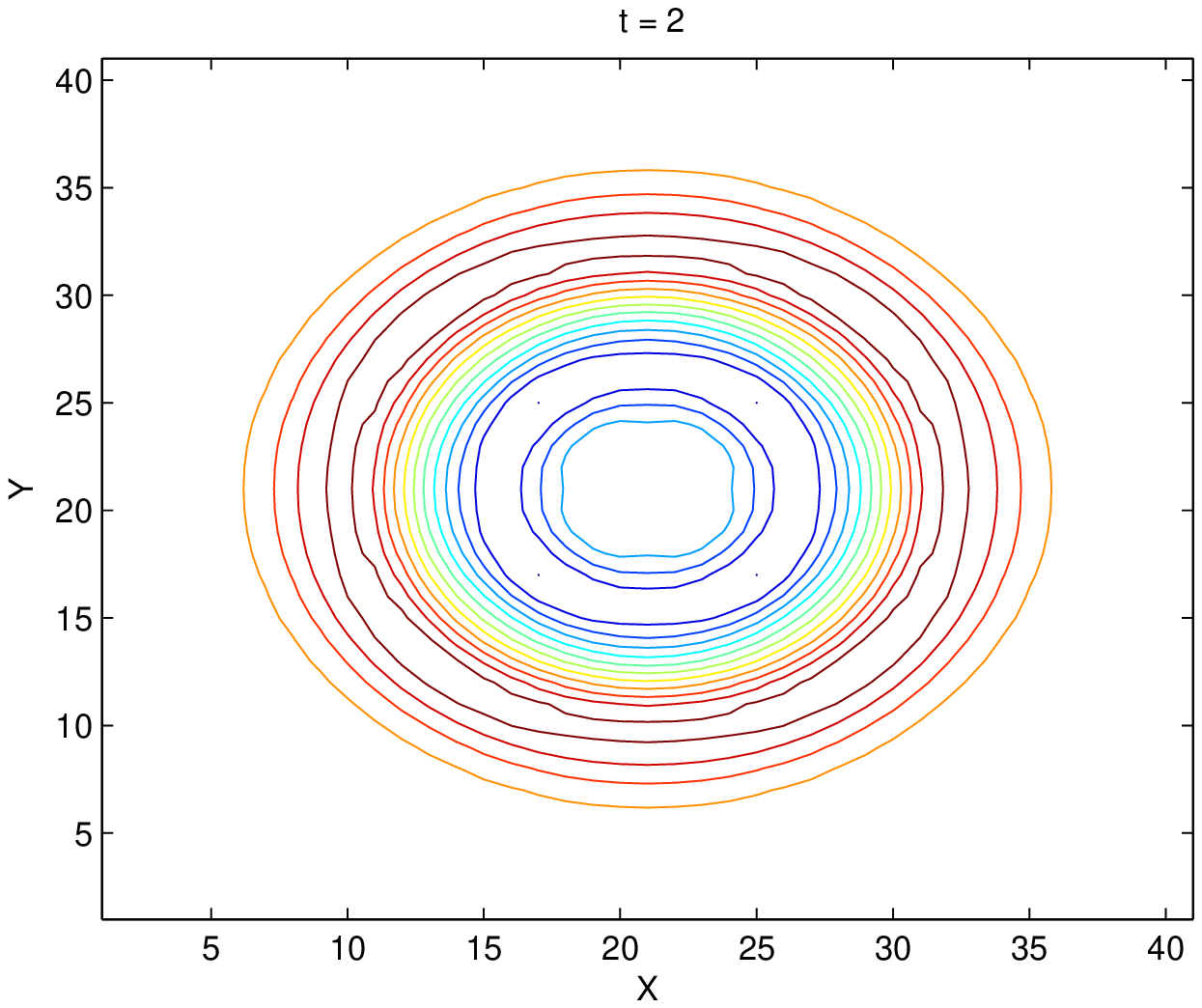}}
\end{subfigure}
\caption{ Numerical solutions of $u$ and the contour profile at $T=2$ with $h=0.1$ and $\Delta t=0.1$.}
\label{Fig78}
\end{figure}

\begin{figure}
\centering
\begin{subfigure}[b]{7cm}
\frame{\includegraphics[width=7cm]{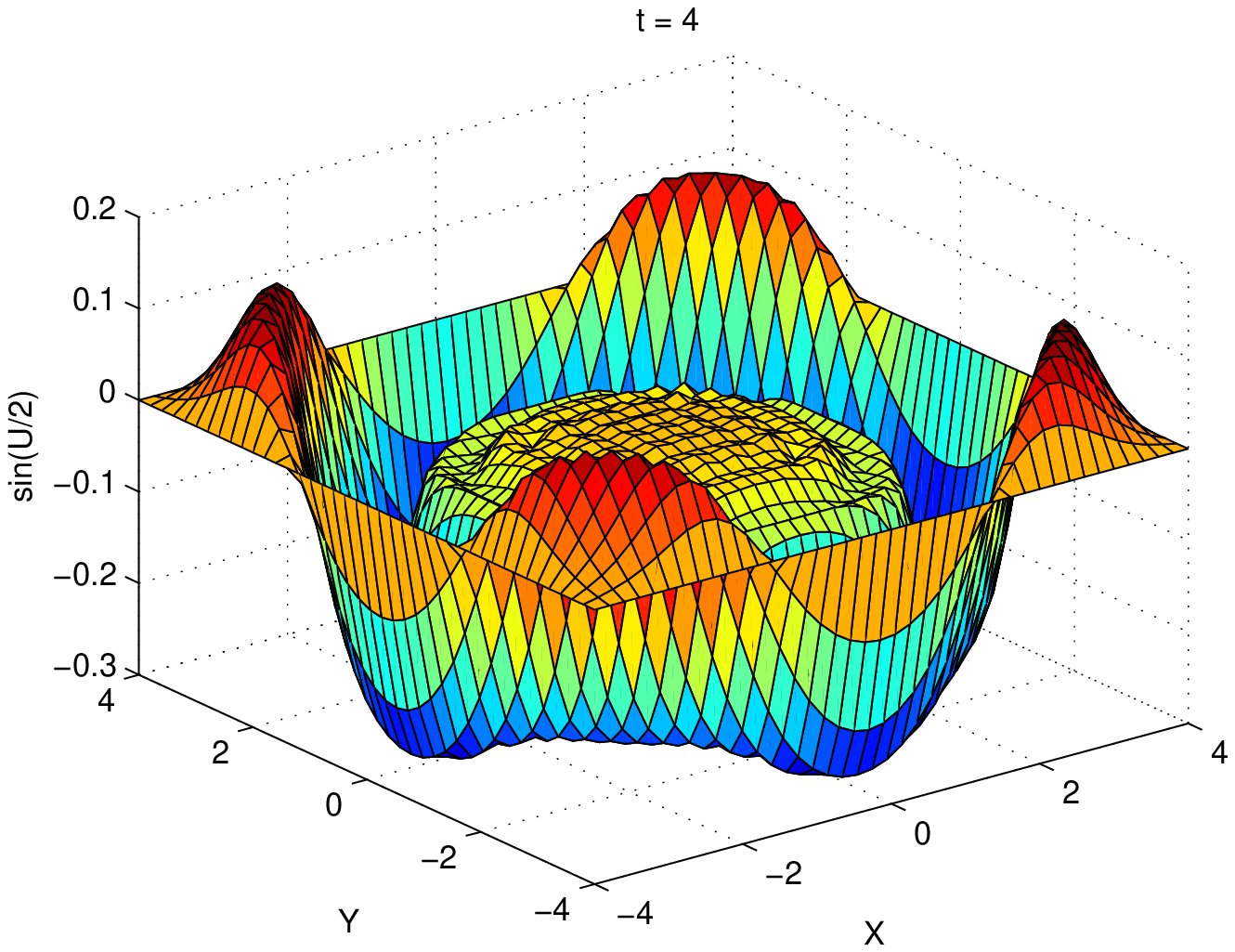}}
\end{subfigure}
\hspace{1cm}
\begin{subfigure}[b]{7cm}
\centering
\frame{\includegraphics[width=7cm]{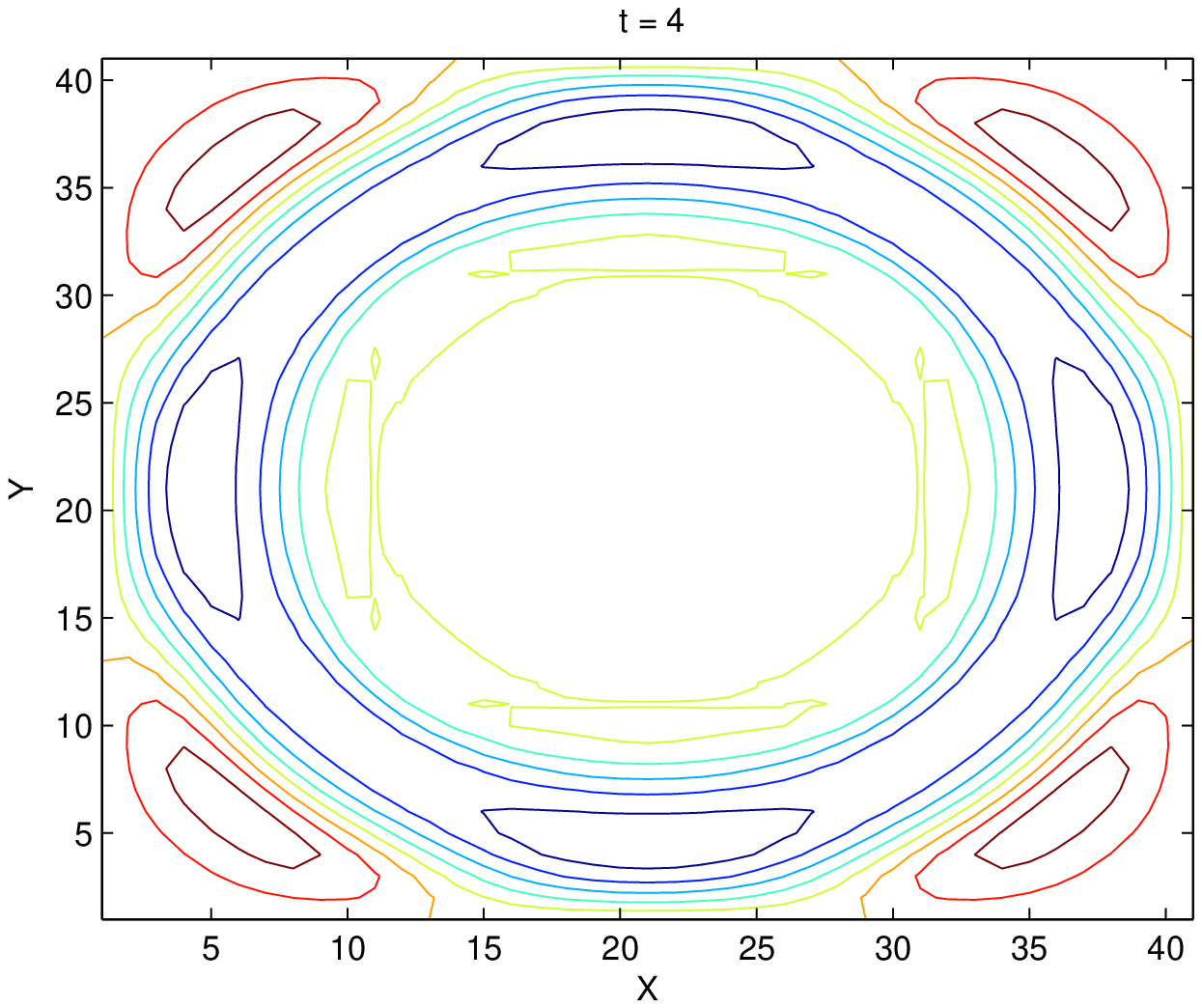}}
\end{subfigure}
\caption{Numerical solutions of $u$ and the contour profile at $T=4$ with $h=0.1$ and $\Delta t=0.1$.}
\label{Fig910}
\end{figure}

\begin{figure}
\centering
\begin{subfigure}[b]{7cm}
\frame{\includegraphics[width=7cm]{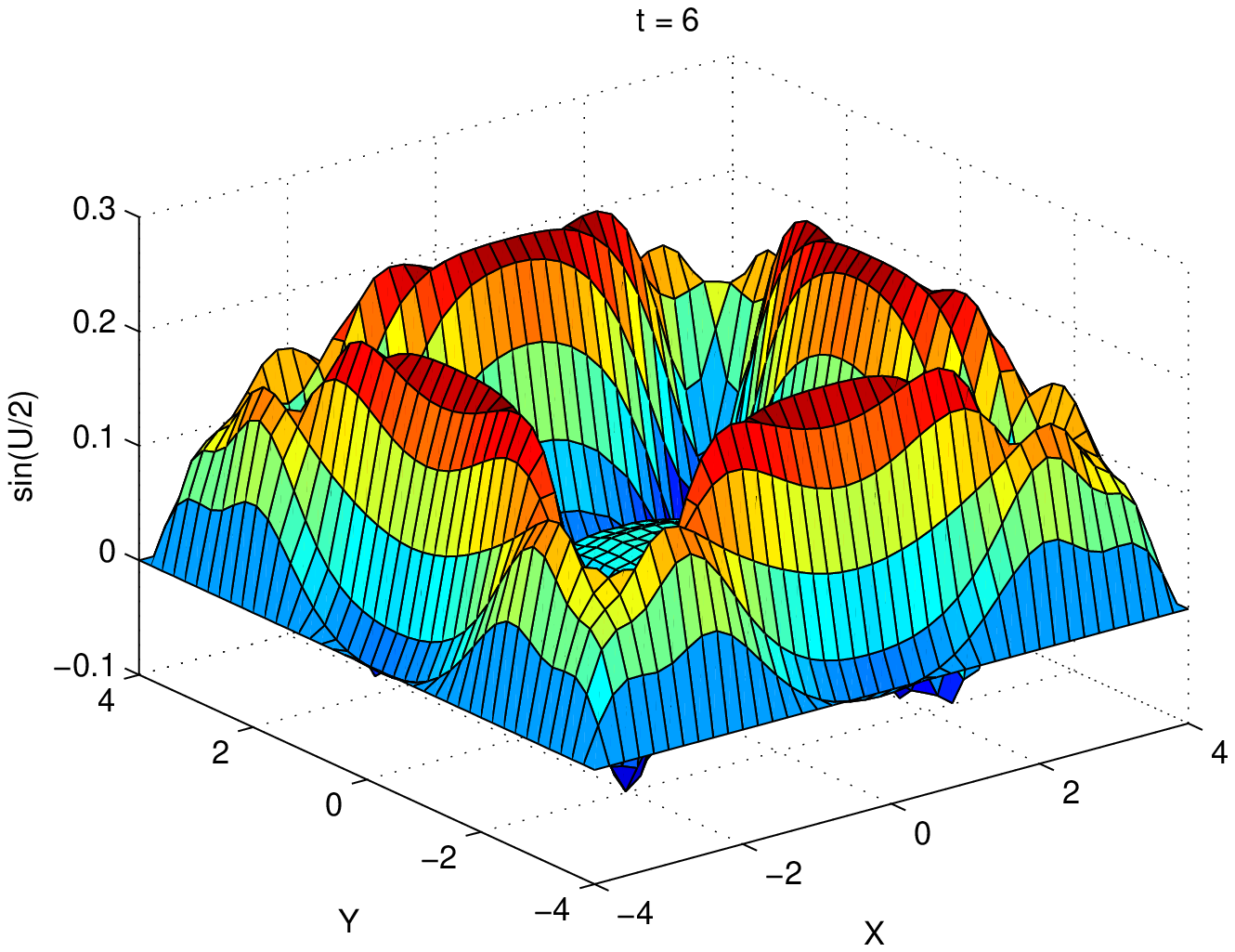}}
\end{subfigure}
\hspace{1cm}
\begin{subfigure}[b]{7cm}
\centering
\frame{\includegraphics[width=7cm]{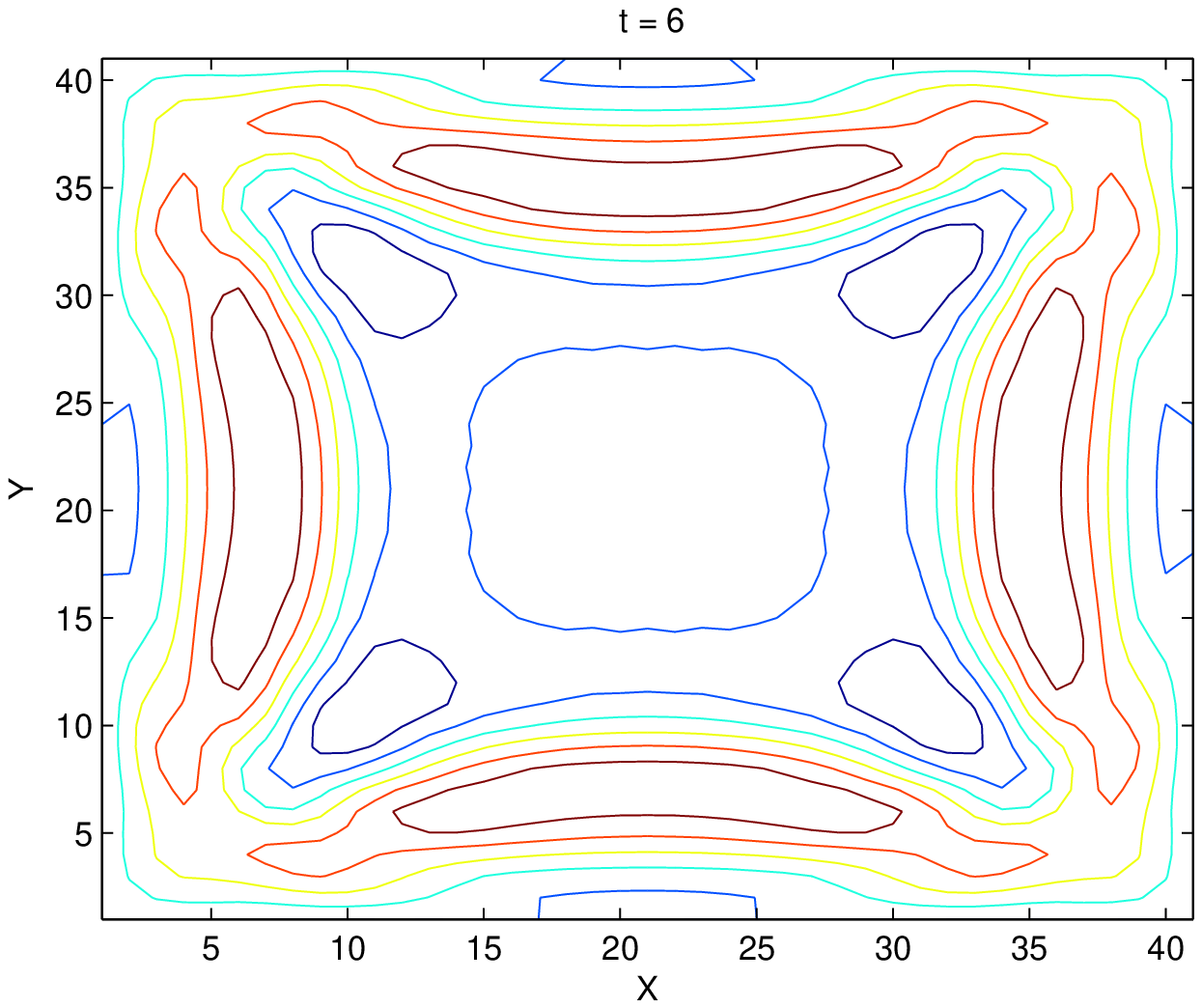}}
\end{subfigure}
\caption{Numerical solutions of $u$ and the contour profile at $T=6$ with $h=0.1$ and $\Delta t=0.1$.}
\label{Fig1112}
\end{figure}

\begin{figure}
\centering
\begin{subfigure}[b]{7cm}
\frame{\includegraphics[width=7cm]{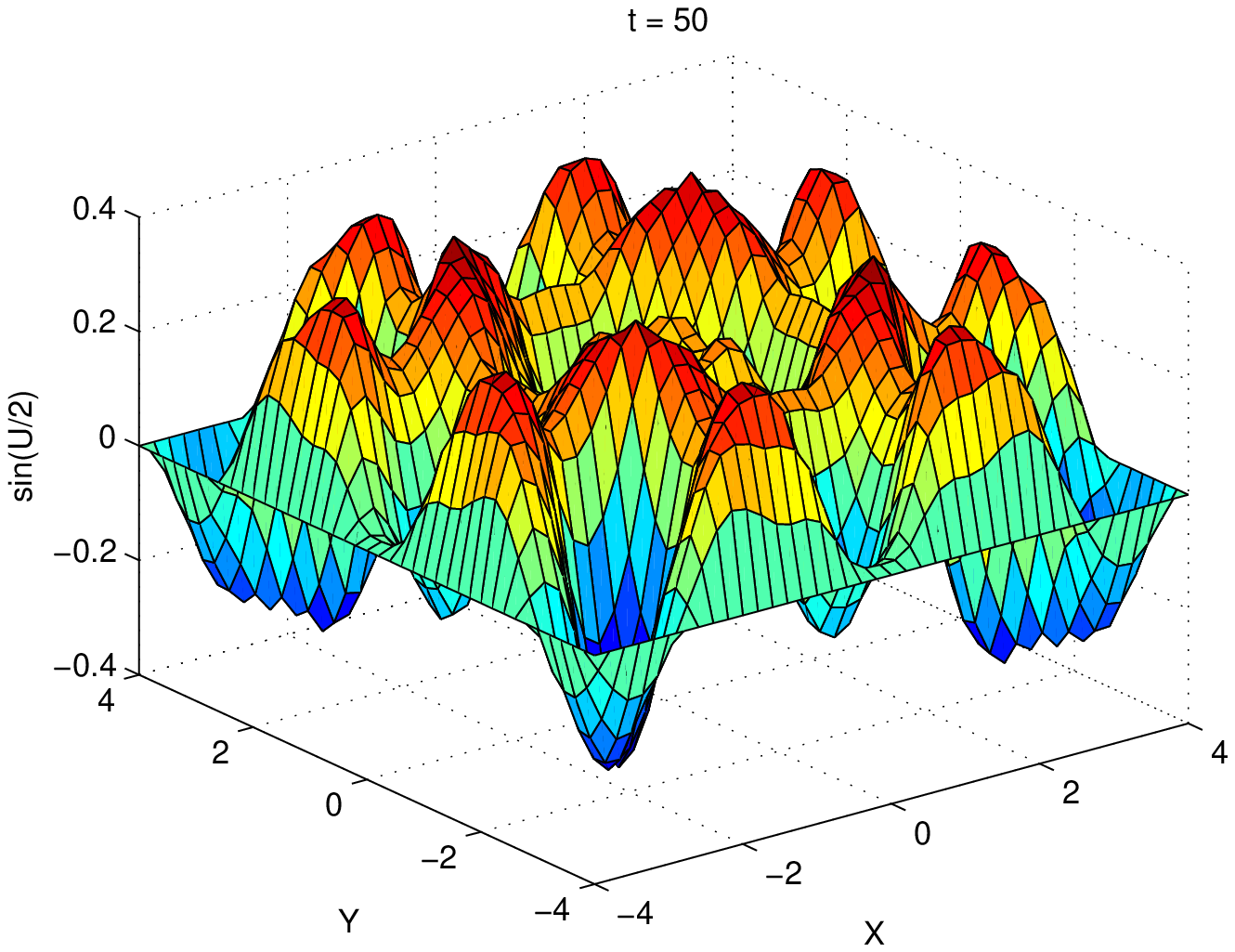}}
\end{subfigure}
\hspace{1cm}
\begin{subfigure}[b]{7cm}
\centering
\frame{\includegraphics[width=7cm]{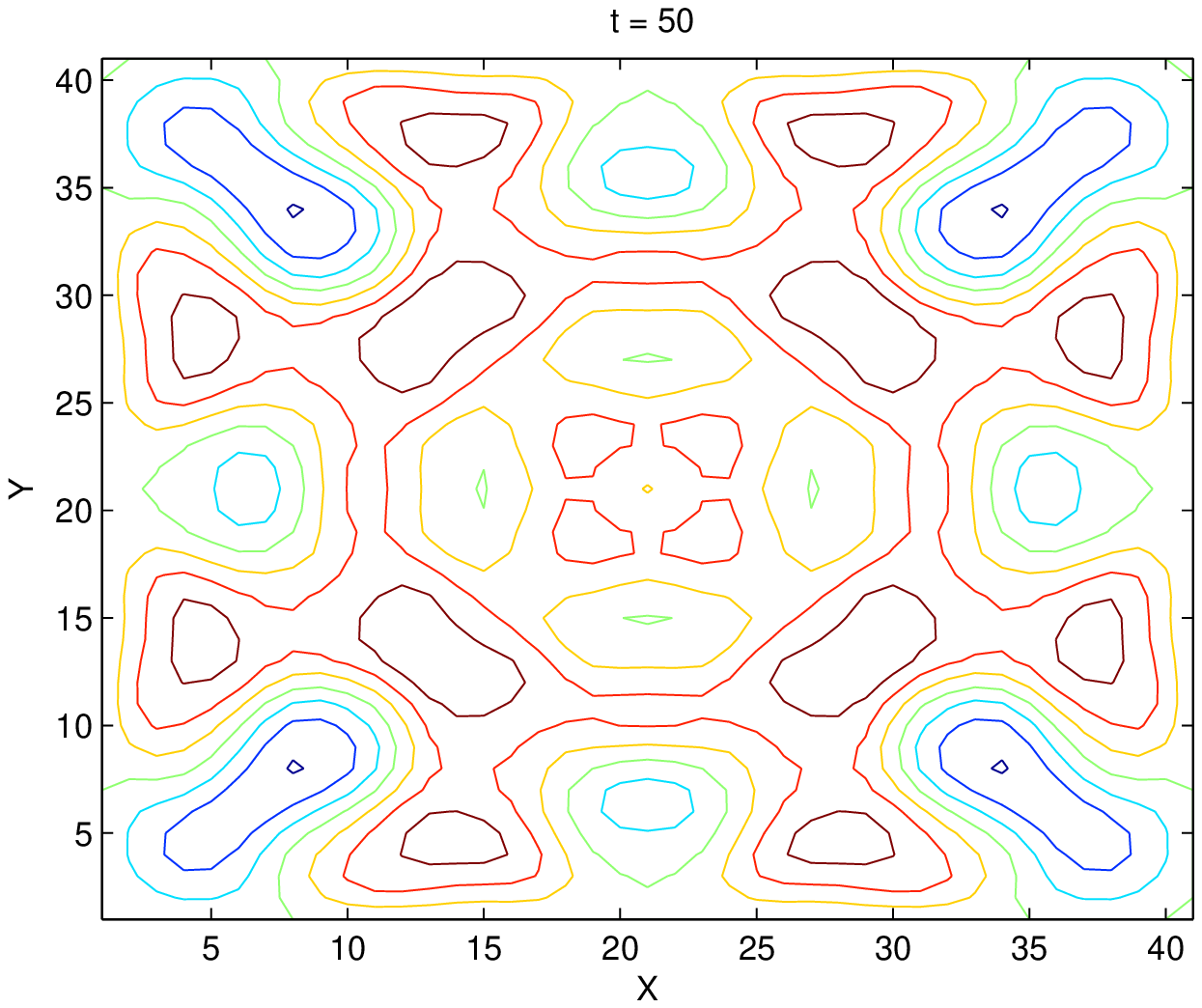}}
\end{subfigure}
\caption{Numerical solutions of $u$ and the contour profile at $T=50$ with $h=0.1$ and $\Delta t=0.1$.}
\label{Fig1314}
\end{figure}

The behavior of a circular ring quasi-soliton arising from the
sine-Gordon equation is named as waves pulsons because of their pulsating behavior. In
this test, we consider the  equation  (\ref{1.1})  on  $ \Omega= [-4, 4] \times [-4, 4]$
with $\phi(x, y) = 1$. The  initial conditions are given by,
\begin{eqnarray}
&&\varphi_1(x,y)=2\arctan(\exp(3-5\sqrt{x^2+y^2})), \label{5.3.1}\\
&&\varphi_2(x,y)=0, \label{5.3.2}
\end{eqnarray}
and the boundary condition is periodic. Similar to \cite{Asg2013}, in order to study the evolution of the ring solitons, we plot both
the surfaces and the corresponding contours in terms of $\sin(\frac{u}{2})$ with $h=0.1$
and $\Delta t=0.1$. As seen from Fig.\ref{Fig56},  Fig.\ref{Fig78}, Fig.\ref{Fig910} and  Fig.\ref{Fig1112},
ring soliton shrinks for initial stage $(t = 0)$, but as time
goes on, oscillations and radiations begin to form and continue to
form up to $t = 4$. At $t = 6$, the graph shows that a ring soliton
is nearly formed again. These graphs are consistent with earlier work on this topic in \cite{Asg2013, Jiw2012}.
Furthermore, with the implicit treatment and the linear iteration algorithm,
it becomes possible to simulate the long time behaviors for such an equation.
In Fig.\ref{Fig1314}, the profile of the numerical solutions of $u$ with  $h=0.1$
and $\Delta t=0.1$ at $T=50$ is presented.

\section{Conclusion}
\label{S6}

In this paper, we discussed a second-order semi-implicit finite difference  scheme for
the 2D sine-Gordon equation, which can admit the discrete energy conservation for the undamped problem.
We also proposed the efficient linear iteration algorithm for approximating the nonlinear
system arising from the implicit treatment of the nonlinear term. Moreover, the iteration
algorithm was proven to be a contraction mapping. In turn, based on truncation errors,
the convergence analysis of the numerical scheme was also shown.
Furthermore, the results of numerical experiments demonstrated
the efficiency and the accuracy of our proposed scheme.

\section*{Acknowledgment}

The authors are very grateful to reviewers for carefully reading this paper and their comments. We appreciate
the support provided for this paper by the Science and Technology Department of Sichuan Province in China (No. 2017GZ0316),
the funds of Sichuan Center for Education Development Research of Education Department (No. CJF15014),
the National Natural Science Funds of China (No. 71471123) and the Fundamental Research Funds for the Central Universities
of China (No. skqy201621).


\end{document}